\documentclass[letterpper,11pt]{article}

\usepackage{amsmath,amsthm}  %
\usepackage{amsfonts}
\usepackage{latexsym}
\usepackage{amssymb}
\usepackage[page]{appendix}

\newtheorem{theorem}{Theorem}

\newtheorem{proposition}[theorem]{Proposition}

\newcommand{\BR}{\bar{\mathbb{R}}}

\DeclareMathOperator{\conv}{conv}

\newtheorem{definition}[theorem]{Definition}

\newcommand{\inner}[2]{\langle #1,#2\rangle}
\newcommand{\norm}[1]{\|{#1}\|}
\newcommand{\R}{\mathbb{R}}

\title{Gauge functions for convex cones}
\author{B. F. Svaiter\footnote{Instituto de Mat\'ematica Pura e Aplicada, 
(IMPA), Estrada Dona
Castorina 110, Rio de Janeiro, RJ, CEP 22460-320, Brazil,
{\tt benar@impa.br}.
The work of this author was partially supported by
CNPq grants no. 
474944/2010-7, 303583/2008-8 and  FAPERJ grant E-26/110.821/2008.
}}

\begin{document}
\maketitle

\begin{abstract}
\noindent
We analyze a class of sublinear functionals which
characterize the interior and the exterior of a convex cone in a normed linear
space.

\medskip

\noindent
keywords: normed vector spaces, cones, sublinear functionals.
\end{abstract}

Let $X$ be a real normed linear space and $K$ a convex closed pointed cone in $X$.
The relation
\begin{equation}
\label{eq:korder}
x\preceq_K y\iff y-x\in K
\end{equation}
is a partial order in $X$ and very often one is interested in
minimizing in the sense of this order a function
$F:\R^n\to X$, that is,
\begin{equation}
\label{eq:p}
\mbox{find } x\in \Omega \mbox{ such that} F(x) \mbox{ is minimal in }F(\Omega),
\end{equation}
where $\Omega$ is a subset of $\R^n$.

Very recently Cauchy method, Newton method and Gradient Projection method
were extended to
problem \eqref{eq:p} in the cases were $F$ is differentiable, the interior of
$K$ is non-empty and
$\Omega=\R^n$ or $\Omega$ is a convex closed set
\cite{MR1778656,MR2108583,MR2049673,MR2515788}.
All these method are $K$-descent method. In particular, the generated sequences
will have $K$-smaller objective values than the initial iterate.
At the hearth of these
  extensions is a \emph{gauge function} for $K$, which
measures how good is a descent direction.
Although the original definition considered cones in finite-dimensional spaces,
its extension to an infinite dimensional setting is straightforward, and this is the case we
shall discuss.
Recall that the \emph{positive polar cone} of $K$ is $K^+$,
\begin{equation}
\label{eq:ppolar}
K^+=\{x^*\in X^*\;|\; \inner{x}{x^*}\geq 0,\;\forall x\in K\}.
\end{equation}

\begin{definition}
\label{def:gf}
A gauge function of/for a closed convex pointed cone $K\subset X$ is
\[
\varphi(x)=\sup_{x^*\in C}\inner{x}{x^*}
\]
where $C\subset K^+\setminus\{0\}$ is a weak-$*$ compact set which generates $K^+$,
in the following sense
\[
K^+=cl \conv\left(\cup_{t\geq 0}tC\right).
\]
\end{definition}

The rationale for the above definition is given in the next elementary result.

\begin{proposition}
\label{pr:bas}
Let $K\subset X$ be a closed convex pointed cone with a nonempty interior. Then
\[
-K=\{x\in X\;|\; \inner{x}{x^*}\leq 0,\; \forall x^*\in K^+\},\quad
(-K)^\circ=\{x\in X\;|\; \inner{x}{x^*}< 0,\; \forall x^*\in K^+\setminus \{0\}\},
\]
where $(-K)^\circ$ stands for the interior of $-K$.
\end{proposition}

\begin{proof}
The inclusion $-K\subset \{x\;|\;\inner{x}{x^*}\leq 0,\;\forall x^*\in K^+\}$ 
follows trivially from \eqref{eq:ppolar}. To show that this inclusion holds as an
equality, suppose that $y\notin -K$. Using Hahn-Banach theorem in its geometric form we conclude that
there exists $x^*$ such that
\[
\inner{x}{x^*}<\inner{y}{x^*},\qquad \forall x\in -K.
\]
As $-K$ is a cone, we conclude that $\inner{x}{x^*}\leq 0$ for any $x\in -K$. Therefore
$x^*\in K^+$ and $y\notin \{x\;|\;\inner{x}{x^*}\leq 0,\;\forall x^*\in K^+\}$.

To prove the second equality, let
\[
V= \{x\;|\; \inner{x}{x^*}< 0,\; \forall x^*\in K^+\setminus \{0\}\}.
\]
If $x\in (-K)^\circ$ and
$\inner{x}{x^*}=0$ for some $x^*\in K^+\setminus \{0\}$, then there exists $y$ in a 
neighborhood of $x$ such that $y\in -K$ and $\inner{y}{x^*}>0$ in contradiction
with the first equality of the proposition.
Therefore, 
\[
(-K)^\circ\subset V\subset -K,
\]
where the second inclusion follows from the definition of $V$ and the first
equality of the proposition.
If the first inclusion is proper, there exists $z\in V\setminus(-K)^\circ$. Using
again Hahn-Banach Theorem we conclude that there exists $x^*$ such that
\[
\inner{x}{x^*}<\inner{z}{x^*}, \qquad \forall x\in(-K)^\circ.
\]
Since the closure of $(-K)^\circ$ is $-K$, $x^*\leq 0$ in $-K$. Therefore $x^*\in K^+\setminus\{0\}$
and, taking the $\sup$ on the left hand-side inequality for $x\in (-K)^\circ$ we conclude that
$\inner{z}{x^*}\geq 0$ in contradiction with the assumption $y\in V$.
\end{proof}

A very natural question is: how general is the class of gauge functions?
Before answering this question in Theorem~\ref{th:gf}, recall that
\emph{Fenchel-Legendre conjugate} of $f:X\to\BR$ is $f^*:X^*\to\BR$,
\[
f^*(x^*)=\sup_{x\in X}\inner{x}{x^*}-f(x)
\]
and the \emph{indicator function} of $A\subset X$ is $\delta_A:X\to\BR$,
\[
\delta_A(x)=
\begin{cases}
0,& x\in A,\\
\infty, &\mbox{otherwise}.
\end{cases}
\]
In view of the above definitions, the function $\varphi$ on Definition~\ref{def:gf} 
can be also expressed as
\begin{equation}
\varphi(x)=(\delta_C)^*(x),\qquad x\in X
\label{eq:gf.alt}
\end{equation}
where we identify $X$ with its canonical injection in to its (topological) bidual.

The aim of this note is to prove that following result:
\begin{theorem}
\label{th:gf}
If $K$ is a closed convex pointed cone with a non-empty interior
then $\varphi$ is a gauge function of/for $K$ if and only if
it satisfies the following properties:
\begin{enumerate}
  \item $\varphi$ is a continuous  sublinear functional
  \item $\varphi(x)<0$ in the interior of $-K$;
  \item $\varphi > 0$ in the complement of $-K$.
\end{enumerate}
\end{theorem}

\begin{proof}
First suppose that $\varphi$ is a gauge function. Then it satisfies trivially item
1. Using Proposition~\ref{pr:bas} (and the weak-$*$ compacity of $C\subset K^+\setminus \{0\}$)
we conclude that $\varphi$ satisfies
item 2. Moreover, since $C$ generate $K^+$, using again Proposition~\ref{pr:bas}
we conclude that
\begin{equation}
  \label{eq:first}
\varphi\geq 0\mbox{ in }X\setminus(-K).
\end{equation}
To show that $\varphi$ satisfies item 3, take $x\in X\setminus(-K)$ and $x_0\in(-K)^\circ$.
Since $-K$ is closed and $x\notin -K$, there exists $\theta\in(0,1)$ such that
\[
x_\theta=(1-\theta)x_0+\theta x\notin -K .
\]
Therefore, using  \eqref{eq:first}, item 1 and item 2 and we conclude that
\[
0\leq\varphi(x_\theta)\leq (1-\theta)\varphi(x_0)+\theta\varphi(x)<\theta\varphi(x),
\]
which proves that $\varphi(x)>0$, that is, $\varphi$ satisfies item 3.

Suppose now that $\varphi$ satisfies items 1, 2 and 3.
Define
\begin{equation}
  \label{eq:C}
C=\{x^*\in X^*\;|\; x^*\leq\varphi\},\qquad
\norm \varphi=\sup_{z\neq 0}\frac{|\varphi(z)|}{\norm{z}}.
\end{equation}
Since $\varphi(x)<0$ for some $x$,
$0\notin C$.
Using items 1 and 2 we conclude that $\varphi\leq 0$ in $-K$. Therefore, if
$x^*\in C$ then $x^*\leq 0\in -K$ and so, $x^*\in K^+$. Altogether we have
\[
C\subset K^+\setminus \{0\}.
\]
If $x^*\in C$ then
$\norm{x^*}=\sup_{z\neq 0}\inner{x}{x^*}/\norm{z}\leq \norm{\varphi}$.
Therefore, $C$ is bounded. 
Since $\varphi$ is sublinear
  \[
  \varphi^*(x^*)=
  \begin{cases}
    0,&x^*\in C\\
    \infty,&\mbox{otherwise}.
  \end{cases}
  \]
  As $\varphi^*$ is convex and weak-$*$ lower semicontinuous, $C$ is convex 
  and weak-$*$ closed.
  Hence $C$ is weak-$*$ compact.
  Moreover, by Fenchel-Moreau Theorem, for any $x\in X$
  \begin{align}
    \label{eq:good}
    \varphi(x)=&\varphi^{**}(x)
    =\sup_{x^*\in C}\inner{x}{x^*}
    =\max_{x^*\in C}\inner{x}{x^*}
  \end{align}
  where the second equality follows from the above expressions for  $\varphi$ and the
  third one from the weak-$*$ compacity of $C$.
  
  Let $V$ be the cone generated by $C$,
  \[
  V=\conv(\cap_{t\geq 0}tC)=
  \cap_{t\geq 0}tC
  \]
  where the second equality follows from the convexity of $C$. As $C\subset K^+$, $V\subset K^+$.
  If $V\neq K^+$, 
  there exists $w_0^*\in K^+\setminus V$ and
  \[
  \{tw_0^*\;|\; t\geq 0\}\cap C=\emptyset.
  \]
  Using the convexity and weak-$*$-compacity of $C$ we conclude that there exists $x_0$ and $b$ such that
  \[
  \inner{x_0}{z^*}<b\leq \inner{x_0}{tw_0^*},\qquad \forall t\geq 0,z^*\in C.
  \]
  Therefore $b\leq 0$, $\inner{x_0}{w_0^*}\geq 0$ and, in view of \eqref{eq:good}, $\varphi(x_0)<0$.
  Hence, by item 2,  $x_0\in (-K)^\circ$. However,
  $w_0^*\in K^+\setminus\{0\}$ and $\inner{x_0}{w_0^*}\geq 0$, in contradiction with Proposition~\ref{pr:bas}.
  So, the assumption $K^+\neq V$ is false and $K^+$ is the cone generated by $C$, which competes the proof that
  $\varphi$ is a gauge function with $C$ given by \eqref{eq:C}.
\end{proof}

In some cases, given a gauge function is easy to provide alternative choices for $C$
(not necessarily weak-$*$ compacts),
as shown in the next example,
Hiriart-Urruty's \emph{oriented distance function}, introduced in \cite{MR543611}
in the framework of nonsmooth scalar optimization in Banach spaces.

\begin{proposition}
  Suppose that $K$ is a pointed convex closed cone with a non-empty interior,
  and let $\varphi:X\to\R$ be given by
  \[
  \varphi(x)=d(x,-K)-d(x,X\setminus -K)
  \]
  Then
  \[
  \varphi(x)=\sup_{x^*\in C}\inner{x}{x^*}
  \]
  for $C=\{x^*\in K^+\;|\; \norm{x^*}=1\}$.
\end{proposition}

\begin{proof}
  Take $x\notin -K$, and let
  $r=d(x,-K)$.
  Define
  \[
  A=-K+B(0,r).
  \]
  Note that $A$ is an open convex set, $0\in A$ and $x\in \bar A$.
  Let $g$ be the Minkowski functional of $A$,
  \[
  g(z)=\inf \{t > 0\;|\; t^{-1}z\in A\},\qquad z\in X.
  \]
  Using Hahn-Banach Theorem we conclude that there exists
 $x^*\in X^*$ such that $x^*\leq g$, $x^*(x)=g(x)=1$.
  Since $-K\subset A$, $x^*\in K^+$. Since $B(0,r)\subset A$,
  $\|x^*\|\leq 1/r$. As $x^*(x)=1$, there exists a sequence $\{y_n\}$ in $-K$
  such that
  \[
  \norm{x-y_n}\to r,\qquad \mbox{ as }n\to\infty
  \]
  Hence
  \[
  \inner{x-y_n}{x^*}=1-\inner{y_n}{x^*}\geq 1.
  \]
  Combining the two above equations we conclude that $\norm{x^*}\geq 1/r$.
  Therefore, $\norm{x^*}=1/r$,
  \[
  rx^*\in C,\qquad \inner{rx^*}{x}=r
  \]
  which proves that
  \[
  \varphi(x)=r\leq \sup_{y^*\in C}\inner{x}{y^*}.
  \]
  To prove the converse, using the definition of $C$ and the sequence
  $\{y_n\}$ we have
  \[
  \sup_{y\in C} \inner{x}{y^*}=
  \sup_{y\in C} \inner{x-y_n}{y^*} + \inner{y_n}{y^*}\leq \norm{x-y_n}
  \]
  so $\varphi(x) = \sup_{y^*\in C}\inner{x}{y^*}$.

  Take $x\in (-K)^\circ$ and let $\rho=d(x,X\setminus -K)$.
  There exists a sequence $\{z_n\}$ in $X\setminus-K$ such that
  \[
  \norm{x-z_n}\to \rho,\qquad \mbox{as }n\to\infty
  \]
  As $z_n\notin -K$, for each $n$ there exists $x^*_n$ such that$\norm{x^*_n}=1$ and
  \[
  \inner{y}{x^*_n}<\inner{z_n}{x^*_n},\quad\forall y\in -K .
  \]
  The above inequality trivially implies $x^*_n\in K^+$, $\inner{z_n}{x^*_n}>0$ and so
  \[
  \inner{x}{x^*_n}=
  \inner{x-z_n}{x^*_n}+\inner{z_n}{x^*_n}
  \geq  -\norm{x-z_n}.
  \]
  Hence
  \[
 \sup_{y^*\in C}\inner{x}{y^*}\geq -\rho=   \varphi(x).
  \]
  To prove the converse, note that for any $x^*\in C$, $\norm{x^*}=1$,
  \[
  \inner{x}{x^*}+\rho=\sup_{\norm{y}<\rho}\inner{x+y}{x^*}\leq\sup_{y'\in -K}\inner{y'}{x^*}= 0
  \]
  which trivially implies the desired inequality.
\end{proof}


\begin{thebibliography}{1}

\bibitem{MR2049673}
L.~M.~Gra{\~n}a Drummond and A.~N. Iusem.
\newblock A projected gradient method for vector optimization problems.
\newblock {\em Comput. Optim. Appl.}, 28(1):5--29, 2004.

\bibitem{MR2515788}
J.~Fliege, L.~M. Gra{\~n}a~Drummond, and B.~F. Svaiter.
\newblock Newton's method for multiobjective optimization.
\newblock {\em SIAM J. Optim.}, 20(2):602--626, 2009.

\bibitem{MR1778656}
J{{\"o}}rg Fliege and Benar~Fux Svaiter.
\newblock Steepest descent methods for multicriteria optimization.
\newblock {\em Math. Methods Oper. Res.}, 51(3):479--494, 2000.

\bibitem{MR2108583}
L.~M. Gra{\~n}a~Drummond and B.~F. Svaiter.
\newblock A steepest descent method for vector optimization.
\newblock {\em J. Comput. Appl. Math.}, 175(2):395--414, 2005.

\bibitem{MR543611}
J.-B. Hiriart-Urruty.
\newblock Tangent cones, generalized gradients and mathematical programming in
  {B}anach spaces.
\newblock {\em Math. Oper. Res.}, 4(1):79--97, 1979.

\end{thebibliography}
\end{document}